\documentclass[11pt,a4paper,reqno]{amsart} 
\usepackage{amsfonts,amsthm, amscd, epsfig, amsmath, amssymb,enumerate}

\numberwithin{equation}{section}

\DeclareMathSymbol{\leqslant}{\mathalpha}{AMSa}{"36} 
\DeclareMathSymbol{\geqslant}{\mathalpha}{AMSa}{"3E} 
\DeclareMathSymbol{\eset}{\mathalpha}{AMSb}{"3F}     
\renewcommand{\leq}{\;\leqslant\;}                   
\renewcommand{\geq}{\;\geqslant\;}                   
\newcommand{\dd}{\,\text{\rm d}}             

\newcommand{\grad}{\nabla} 

\newcommand{\la}{\label} 
\newcommand{\eps}{\epsilon} 
\newcommand{\be}{\begin{equation}}

\def\1{\ifmmode {1\hskip -3pt \rm{I}} \else {\hbox {$1\hskip -3pt \rm{I}$}}\fi}

\newtheorem{Th}{Theorem}[section] 
\newtheorem{Le}[Th]{Lemma} 
\newtheorem{Pro}[Th]{Proposition}

\newcommand{\cD}{\ensuremath{\mathcal D}} 
\newcommand{\cE}{\ensuremath{\mathcal E}} 
\newcommand{\cF}{\ensuremath{\mathcal F}} 
\newcommand{\cG}{\ensuremath{\mathcal G}}

\newcommand{\cL}{\ensuremath{\mathcal L}}

\newcommand{\bbN}{{\ensuremath{\mathbb N}} }

\newcommand{\bbR}{{\ensuremath{\mathbb R}} }

\newcommand{\si}{\sigma} 
\newcommand{\ent}{{\rm Ent} } 
 
\newcommand{\dert}{\frac{{\rm d}}{{\rm d} t}}

\let\a=\alpha    \let\d=\delta  
 \let\g=\gamma       \let\l=\lambda
            
\let\r=\rho      
  
   \let\G=\Gamma   
\let\O=\Omega      

\def\\{\hfill\break}

\def\thsp{\thinspace}

\def\tthsp{\kern .083333 em}

\def\?{\mskip -10mu}
\def\indbox#1{\hbox to \parindent{\hfil\ #1\hfil} }

\newfam\msafam
\newfam\msbfam
\newfam\eufmfam
\def\hexnumber#1{%
  \ifcase#1 0\or 1\or 2\or 3\or 4\or 5\or 6\or 7\or 8\or
  9\or A\or B\or C\or D\or E\or F\fi}
\font\tenmsa=msam10
\font\sevenmsa=msam7
\font\fivemsa=msam5
\textfont\msafam=\tenmsa
\scriptfont\msafam=\sevenmsa
\scriptscriptfont\msafam=\fivemsa        
\edef\msafamhexnumber{\hexnumber\msafam}%
\mathchardef\restriction"1\msafamhexnumber16
\mathchardef\ssim"0218
\mathchardef\square"0\msafamhexnumber03
\mathchardef\eqd"3\msafamhexnumber2C
\def\QED{\ifhmode\unskip\nobreak\fi\quad
  \ifmmode\square\else$\square$\fi}            
\font\tenmsb=msbm10
\font\sevenmsb=msbm7
\font\fivemsb=msbm5
\textfont\msbfam=\tenmsb
\scriptfont\msbfam=\sevenmsb
\scriptscriptfont\msbfam=\fivemsb

\font\teneufm=eufm10
\font\seveneufm=eufm7
\font\fiveeufm=eufm5
\textfont\eufmfam=\teneufm
\scriptfont\eufmfam=\seveneufm
\scriptscriptfont\eufmfam=\fiveeufm

\def\({\left(}
\def\){\right)}
\let\neper=e
\let\ii=i

\def\nep#1{ \neper^{#1}}

\def\tc{\thsp | \thsp}
\let\<=\langle
\let\>=\rangle

\outer\def\nproclaim#1 [#2]#3. #4\par{\medbreak \noindent
   \talato(#2){\bf #1 \Thm[#2]#3.\enspace }%
   {\sl #4\par }\ifdim \lastskip <\medskipamount 
   \removelastskip \penalty 55\medskip \fi}
\def\thmm[#1]{#1}
\def\teo[#1]{#1}
\def\sttilde#1{%
\dimen2=\fontdimen5\textfont0
\setbox0=\hbox{$\mathchar"7E$}
\setbox1=\hbox{$\scriptstyle #1$}
\dimen0=\wd0
\dimen1=\wd1
\advance\dimen1 by -\dimen0
\divide\dimen1 by 2
\vbox{\offinterlineskip%
   \moveright\dimen1 \box0 \kern - \dimen2\box1}
}
\begin{document}

\title[Entropy dissipation estimates in a zero--range dynamics]
{Entropy dissipation estimates in a zero--range dynamics}

\begin{abstract}
We study the exponential decay of relative
entropy functionals for zero--range processes
on the complete graph. 
For the standard model with rates increasing at infinity
we prove entropy dissipation estimates,
uniformly over the number of particles and the number of vertices.

\bigskip



\end{abstract}

\author[P. Caputo]{Pietro Caputo}
\address{Dip. Matematica, Universita' di Roma Tre, L.go S. Murialdo 1,
00146 Roma, Italy} \email{caputo\@@mat.uniroma3.it}
\author[G. Posta]{Gustavo Posta} \address{Dip. Matematica, 
Politecnico di Milano, 
P.za Leonardo da Vinci 32, 
I-20133 Milano, Italy} 
\email{gustavo.posta\@@polimi.it}



\maketitle

\thispagestyle{empty}

\section{Introduction, Models and Results}
Functional estimates such as Poincar\'e and logarithmic Sobolev
inequalities play an important role 
in the study of approach to stationarity for Markov semigroups,
see e.g.\ \cite{Mar} for a recent survey.
Logarithmic Sobolev inequalities are known to imply 
exponential decay of relative entropy which in turn provides a natural 
way to bound mixing times in total variation norm \cite{DSC}.
%
As we can see already in simple birth--and--death processes,
however, in discrete settings
logarithmic Sobolev inequalities may become an unnecessarily strong
requirement if we are interested in decay to equilibrium
in relative entropy or total variation. 
Motivated by this observation, modified versions 
of the logarithmic Sobolev inequality have been recently proposed 
and studied by several authors \cite{BobLed,DPP,GQ,Goel,BobTet}. 
As emphasized in \cite{DPP,GQ,Goel,BobTet} a key estimate is the one 
relating directly the relative entropy functional and its time--derivative
along the semigroup.
Such entropy dissipation inequalities have been extensively studied 
in the literature on the approach to 
equilibrium for the Boltzmann equation, see \cite{Vil}
and references therein. 
Our aim in this paper is to investigate 
the validity of entropy dissipation bounds for some  
models of interacting random walks on the complete graph
known as zero--range processes.

\smallskip

The complete graph zero--range dynamics is the 
continuous time Markov chain
described as follows. For each positive integer $L$ we consider the
set of vertices $V_L=\{1,\dots,L\}$,
the state space is the product $\O_L=\bbN^L$
and a configuration $\eta\in\O_L$ is interpreted as an occupation number
vector, i.e.\ $\eta_x$ is the number of particles at $x\in V_L$.
At each vertex $x\in V_L$ we associate 
a rate function $c_x:\bbN\to \bbR$ such that $c_x(0)=0$ and $c_x(n)>0$
for every $n\geq 1$. We often extend $c_x$ to a function on $\O_L$ 
by setting $c_x(\eta)=c_x(\eta_x)$. 
Every vertex $x\in V_L$ waits an exponentially distributed time
with mean $1/c_x$ before one particle is moved from $x$ to a 
uniformly chosen vertex of $V_L$. More precisely, the 
Markov generator is given by
\be
\cL f=\frac1L\sum_{x,y} c_x\grad_{xy}f\,,
\la{gen}
\end{equation}
with the sum extending over all $x,y\in V_L$. Here 
$\grad_{xy}f$ stands for the gradient $f^{xy} -f$, 
with $f^{xy}(\eta) = f(\eta^{xy})$, $\eta^{xy}$ being 
the configuration in which a 
particle has been moved from $x$ to $y$, i.e.\ 
$(\eta^{xy})_x = \eta_x-1$, $(\eta^{xy})_y = \eta_y+1$, and
$(\eta^{xy})_z = \eta_z$, $z\neq x,y$.
We agree that $\eta^{xy}=\eta$, if $\eta_x=0$.  
Note that if the functions $c_x$ were all 
linear, i.e. $c_x(n)=\l_x n$, $\l_x>0$, 
the resulting random walks on the complete
graph with $L$ vertices would be independent. 
The interaction is therefore hidden in the non--linearity of $c_x$
and has zero--range in the sense that jump rates out
of $x$ are only determined by
the configuration at $x$. 
The process is reversible w.r.t.\ the product measure 
$\mu_L(\eta)=\prod_{x\in V_L}\mu_x$, where 
$\mu_x$ is the probability on $\bbN$ 
given by 
\be
\mu_x(0)=\frac1{Z_x}\,,\quad \mu_x(n)=\frac1{Z_x}\prod_{k=1}^n 
\frac1{c_x(k)}\,.
\la{mux}
\end{equation}
Since the process conserves the initial number of particles,
letting $\nu:=\nu_{L,N}$ 
denote the probability $\mu_L$ conditioned on the event 
$N=\sum_{x\in V_L}\eta_x$, we obtain, for every $N\geq 1$ and $L\geq 2$, an
irreducible finite state Markov chain with reversible measure $\nu_{L,N}$. The associated 
Dirichlet form is given by 
\be
\cE_\nu(f,g) = - \nu\left[f (\cL g)\right]
= \frac1{2L} \sum_{x,y} \nu\left[c_x \grad_{xy}f\grad_{xy}g
\right]\,,
\la{Dirich}
\end{equation}
where $f,g$ are arbitrary functions and 
the notation $\nu[f]$ is used for the 
expectation $\int f\mathrm{d}\nu$.  
Local variants of the zero--range dynamics
have been considered in the literature, 
especially in connection with  
hydrodynamical limits \cite{KL}. 
If we allow, for instance, a
particle at $x$ to jump to $x+1$ or $x-1$ only, we have the local 
Dirichlet form
\be
\cD_\nu(f,g) 
= \frac1{2} \sum_{x=1}^{L-1} \nu\left[c_x \grad_{x,x+1}f\grad_{x,x+1}g
\right]\,.
\la{Dirich2}
\end{equation}
Because of the permutation symmetry
of the model it is natural to study the complete graph
dynamics, which is more tractable from
the analytical point of view. 
Moreover, it turns out that in some cases
sharp estimates on the decay to equilibrium for the local variants 
are deduced from the corresponding bounds on the complete graph,
see e.g.\ \cite{C,Morris}.

\smallskip

Let us now recall the notion of entropy and the associated inequalities.
%
As usual the 
the entropy of a function $f\geq 0$ is written $\ent_\nu(f)=
\nu[f\log f] - \nu[f]\log\nu[f]$.
When 
$f\geq 0$ and
$\nu[f]=1$, $\ent_\nu(f)$ coincides with the relative entropy of the
probability $\nu f$ 
w.r.t.\ $\nu$.
%
%
%
Setting $f_t = \nep{t\cL}f$ we have 
\be
\frac{\dd}{\dd t} \,\ent_{\nu}(f_t) = -\, \cE_\nu(f_t,\log f_t)\,.
\la{entd1}
\end{equation}
Therefore the entropy dissipation constant
\be
\g(L,N)=\sup_{f > 0} \,\frac{\ent_\nu(f)}{\cE_\nu(f,\log f)}\,,
\la{entd2}
\end{equation}
is the best constant $\g$ such that 
\be
\ent_\nu(f_t)\leq \nep{-t/\g} \,\ent_\nu(f)\,,
\la{deco}
\end{equation}
for every non--negative function $f$. 

Since (see e.g.\ \cite{BobTet,Goel})
\be
\cE_\nu(f,\log f)\geq 4\,\cE_\nu(\sqrt f,\sqrt f)\,,
\la{decos}
\end{equation}
we see that (\ref{deco}) is implied by the usual logarithmic Sobolev
inequality. Namely, if $s(L,N)$ denotes the logarithmic Sobolev
constant defined by (\ref{entd2}) with $\cE_\nu(f,\log f)$
replaced by $\cE_\nu(\sqrt f,\sqrt f)$, then $\g(L,N)\leq s(L,N)/4$.  
The name 
``modified'' logarithmic Sobolev
constant is sometimes used for $\g(L,N)$. 
%
Note that $s(L,N)$ can be much larger than $\g(L,N)$. 
As an example, consider the simple random walk on the complete
graph with $L$ vertices, which corresponds to the case $N=1$ with
homogeneous rates: $c_x = c_y$, all $x,y\in V_L$. 
Simple computations show that in this case
the logarithmic Sobolev constant $s(L,1)$
grows with $L$ as $\log L$ while $\g(L,1)$ remains bounded.
%

\smallskip

Our main result is obtained under the hypothesis of 
homogeneous Lipschitz rates increasing at infinity.
We formulate this as follows.

\smallskip

There exists $c:\bbN\to \bbR_+$ such that $c_x(n)=c_y(n)=c(n)$, for all
$x,y$ and $n\in\bbN$. Moreover $c(0)=0$, $c(k)>0$ for every
$k\geq1$, and there exist
$C<\infty$, $\d>0$ and $n_0\in\bbN$ such that 
$$
c(m)-c(n)\geq \d\,,
\eqno({\rm H}.1)
$$
for every $n\in\bbN$ and $m\geq n+n_0$, and  
$$
\sup_{n\geq 0} |c(n+1)-c(n)|\leq C\,.
\eqno({\rm H}.2)
$$

\smallskip

Zero--range processes satisfying (H.1) and (H.2) 
have been extensively studied \cite{LSV,C,DP,JLQY}.
Using a version of the Lu--Yau martingale approach \cite{LY}, 
Landim, Sethuraman and Varadhan \cite{LSV} proved that the spectral gap
of the local zero--range dynamics scales diffusively. 
For the complete graph model a uniform  
spectral gap estimate was proved in \cite{C}
following the Carlen--Carvalho--Loss approach to the determination of the 
spectral gap \cite{CCL}.
Using a version of the Cancrini--Martinelli duplication method \cite{CM,CMR}
Dai Pra and Posta \cite{DP1,DP}
have recently established a diffusive estimate for the 
logarithmic Sobolev constant of the local dynamics.

We will show 
that under (H.1) and (H.2) 
the entropy dissipation constant $\g(L,N)$ has a uniform upper bound.
\begin{Th}
\la{mainre}
Assume (H.1) and (H.2). Then 
\be
\sup_{L\geq 2} \,\sup_{N\geq 1}\, \g(L,N) \,<\,\infty\,.
\la{Ug}
\end{equation}
\end{Th}
We conclude with some remarks on this result and on the organization of
the paper.

\smallskip\smallskip

\noindent 
1. As in \cite{LSV,C,DP} our estimate is 
uniform over the number of particles $N$. 
This uniformity 
can no longer be expected if one drops the assumption (H.1),
see for instance \cite{Morris} where the spectral gap of the 
complete graph model with constant rates is shown to be of order
$L^2/(N^2+L^2)$. 

\smallskip\smallskip

\noindent 
2. Following the standard martingale approach, our proof consists
in setting up a recursion on the number of vertices $L$. The first
step requires a one--vertex entropy dissipation 
estimate. This is established in section 2 as a consequence of 
a more general one--dimensional bound for log--concave measures on
$\bbN$ which is of independent interest. 
The rest of the proof is given in section 3. Here we need to
adapt techniques developed for spectral gap and logarithmic
Sobolev inequalities to the more delicate entropy dissipation
estimate. In particular, the one--vertex bound is shown to
produce certain covariance terms in subsection 3.1. The crucial bound
on these covariances is established  
in subsections 3.2--3.4 by combining  the methods of \cite{LPY} and
\cite{DP}. 

\smallskip\smallskip

\noindent 
3. It is natural to try to extend the result to the case of
   inhomogeneous rates. 
One can consider the rates $c_x(n)=\l_x c(n)$
with constants $\l_x\in[a_1,a_2]$ for some $0<a_1<a_2<\infty$
and $c(\cdot)$ satisfying (H.1) and (H.2). For this model
we conjecture that there exists $C<\infty$ depending only on
$a_1,a_2$ and the constants appearing in (H.1), (H.2) such that 
\be
\la{glnc}
\g(L,N)\leq C\,,
\end{equation}
uniformly in $N,L$.  
   We see no serious 
difficulty in adapting the one--vertex estimate of section 2
   and the arguments up to and including subsection
   3.1 to this more general setting. On the other hand,  
the extension of the covariance estimate to this case seems to be more
delicate.  In analogy with \cite{LN} one can take the $\bar\l=\{\l_x\}$ as the
realization of an i.i.d.\ random environment. The covariance terms produce
new terms involving fluctuations of the environment. Combining the
   strategy of \cite{LN} with our arguments in sesction 3
one could possibly
show that the bound (\ref{glnc}) holds with a constant $C= C(\bar \l)$ depending on
the random field and such that $C <\infty$ almost surely. 
Establishing the full conjecture (\ref{glnc}) however  
seems to be a more challenging problem which deserves further investigation. 
If the rates are assumed to be pointwise increasing
a first result is the following
perturbative bound we obtained in \cite{abcdp} by a discrete
version of the Bakry--Emery $\G_2$ criterium. Suppose the rates $c_x$ are
arbitrary functions satisfying:
there exist $\l>\d>0$ such that 
\be
\l\leq c_x(n+1)-c_x(n)\leq \l+\d\,,
\la{inh}
\end{equation}
for every $x$ and $n$;
then $\g(L,N)\leq (\l-\d)^{-1}$ for all $L\geq 2$ and $N\geq 1$.
We believe that a uniform estimate on $\g(L,N)$ as in (\ref{Ug})
should hold under (\ref{inh}) for any $\d>0$
without the restriction $\d<\l$. However, as explained in \cite{abcdp} the $\G_2$
approach breaks down when there is no restriction on $\d$.
For the spectral gap the situation is easier.
In fact, as shown in \cite{abcdp},  assuming (\ref{inh}) one has that 
the spectral gap is bounded below by $\l$ independently of $\d$.

\section{One--vertex estimates}
The goal of this section is to show that
the one--vertex marginal of the canonical measure 
$\nu$ satisfies a uniform entropy dissipation bound, see Proposition 
\ref{one-v} below. From now on the rate function $c:\bbN\to\bbR$ is assumed
to satisfy the conditions (H.1) and (H.2). For any $L\geq 2$ and $N\geq 1$
we write as usual $\nu= \nu_{L,N}$  for the 
homogeneous zero--range canonical measure associated to the rate function 
$c$. We also write $\nu_x$ for the marginal of $\nu$ at $x$, i.e.\
$\nu_x(n) = \nu(\eta_x=n)$. 
\begin{Pro}
\la{one-v}
There exists $C<\infty$ such that,  
for any $L\geq 2, N\geq 1$, $x\in V_L$
and for any function $u:\bbN\to\bbR_+$
with $\nu_x[u]=1$ we have
\be
\sum_{n=0}^N\nu_x(n) u(n)\log u(n) \leq 
C\, \sum_{n=0}^N
\nu_x(n) c(n) [u(n)- u(n-1)] \log \frac{u(n)}{u(n-1)}\,.
\la{22}
\end{equation}
\end{Pro}
The proof will be 
based on the one--dimensional estimate established in \cite{abcdp}, 
which we recall below.
\subsection{A one--dimensional estimate}
Let
$\mu:\bbN\to[0,1]$ be a probability vector and consider the 
birth and death process with birth rate $r_+(n)$ and death rate
$r_-(n)$ satisfying the detailed balance w.r.t.\ $\mu$:
$$
r_-(n)\mu(n) = r_+(n-1)\mu(n-1)\,,\quad n\geq 1\,.
$$
We assume $r_-(0)=0$. 
The following estimate can be found in \cite{abcdp}.
\begin{Le}
\la{claim1}
Let $r_-,r_+$ satisfy 
\begin{gather}
r_-(n+1)-r_-(n)\geq \d_-\\
r_+(n)-r_+(n+1)\geq \d_+\,,
\end{gather}
with some constants $\d_-,\d_+\geq0$. Then, for every $u\geq 0$ such that $\mu[u]=1$
we have 
\be
\sum_{n=0}^\infty \mu(n)\,u(n)\,\log u(n) \leq \d^{-1}
\sum_{n=0}^\infty \mu(n)\,r_-(n)\,[u(n)- u(n-1)]  
\log \frac{u(n)}{u(n-1)}\,,
\la{entdiss}
\end{equation}
where $\d:=\d_-+\d_+$. 
\end{Le}

\subsection{1D--equivalence with the case of increasing rates}
The next step is the following equivalence lemma, whose 
proof can be found in \cite{abcdp}. 
Let $n_0$ be the constant appearing in (H.1). We 
define
\be
\tilde c(k) = c(k) + \frac1{n_0}\sum_{j=1}^{n_0-1}
\frac{n_0-j}{n_0}\,[c(k+j) +c(k-j)- 2c(k)]\,,\quad \;k\geq n_0\,.
\la{ctilde}
\end{equation}
When $k<n_0$ we simply set $\tilde c(k)= \tilde c(n_0)k/{n_0}$.
Let us call $\tilde \mu$ the one--coordinate 
zero--range measure obtained from $\tilde c$, i.e.
\be
\tilde \mu(0)=\frac1{\tilde Z}\,,\quad \tilde \mu(n)=
\frac1{\tilde Z}\prod_{k=1}^n 
\frac1{\tilde c(k)}\,.
\la{mutilde}
\end{equation}
\begin{Le} 
\la{equi}
The rate function $\tilde c$ is uniformly increasing:
there exists $\d>0$ such that  for every $k\in\bbN$
\be
\tilde c(k+1) - \tilde c(k) \geq \d\,.
\la{equi0}
\end{equation}
Moreover, $\mu$ and $\tilde\mu$ are equivalent: 
there exists $C\in[1,\infty)$ such that for every $n\in\bbN$
\be
\frac1C \,\leq\, \frac{\tilde \mu(n)}{\mu(n)}\, \leq\, C\,.
\la{equi1}
\end{equation}
\end{Le}

\subsection{Proof of Proposition \ref{one-v}}
From a standard comparison result (see e.g.\
\cite{Led}, Lemma 1.2) Proposition \ref{one-v} follows 
if we can prove that $\nu_x$ is equivalent to a 
probability $\hat\nu_x$ on $\bbN$ for which 
the estimate (\ref{22}) is known to hold. Here equivalence 
means a double bound as in (\ref{equi1}). Since this notion will be used 
repeatedly in what follows we introduce a special notation for it:
We say that $a:\bbN\to\bbR_+$ is equivalent to $b:\bbN\to\bbR_+$ 
and write $a\asymp b$ whenever 
there exists a universal constant $C\in[1,\infty)$
(independent of $L$ and $N$) such that 
$C^{-1}\leq a/b \leq C$.

Recall the notation $\mu_L$
for the product $\otimes_{x\in V_L}\mu_x$. We shall use the shortcut notation
$\mu_L(k)$ for the probability of the event $\sum_{j=1}^{L}\eta_j=k$,
for every $L\geq 2$ and $k\in\bbN$. By definition
\be
\nu_x(n) = 
\mu_x(n)\,
\frac{\mu_{L-1}\left(N-n\right)}
{\mu_{L}\left(N\right)}\,.
\la{nuxn}
\end{equation}
Let $\tilde \mu_x$ denote the one--vertex measure with rate $\tilde c$
given by (\ref{ctilde}) and write 
$\tilde \mu_L = \tilde \mu_x \otimes (\otimes_{y\in V_L\setminus \{x\}}\mu_y)$.
From Lemma \ref{equi} we know that 
$\mu_x\asymp \tilde \mu_x$ and $\mu_L\asymp \tilde \mu_L$.
Therefore $\nu_x\asymp \tilde \nu_x$ where 
\be
\tilde \nu_x(n) = 
\tilde \mu_x(n)\,
\frac{\mu_{L-1}\left(N-n\right)}
{\tilde \mu_{L}\left(N\right)}\,.
\la{nutx}
\end{equation}
We will use the following lemma.
\begin{Le}\la{lele}
Let $\hat\nu_x$ be a probability on $\{0,1,\dots,N\}$ such that
the function 
\be
V(n):= - \log \,\frac{\hat\nu_x(n)}{\tilde \mu_x(n)} 
\la{vn}
\end{equation}
satisfies
\be
\grad^2 V(n) = V(n+2) + V(n) - 2 V(n+1) \geq 0\,,\quad\; n=0,1,\dots,N-2\,.
\la{logco}
\end{equation}
Then, for every function $u:\bbN\to\bbR_+$
with $\hat \nu_x[u]=1$ we have
\be
\sum_{n=0}^N\hat \nu_x(n) u(n)\log u(n) \leq 
C\, \sum_{n=1}^N
\hat \nu_x(n) c(n) [u(n)- u(n-1)] \log \frac{u(n)}{u(n-1)}\,.
\la{222}
\end{equation}
where $C$ is a constant depending only on the parameters
appearing in (H.1) and (H.2).
\end{Le}
\proof
We extend $\hat \nu_x$ to a probability on $\bbN$ by
setting $\hat\nu_x(k)=0$, $k\geq N+1$. 
We apply Lemma \ref{claim1} with $\mu=\hat\nu_x$,
$r_-(n)=\tilde c(n)$. Then, by reversibility and (\ref{vn}):
$$
r_+(n) = \tilde c(n+1) \,\frac{\hat\nu_x(n+1)}{\hat\nu_x(n)}
= \nep{ - \grad V(n)}\,.
$$
By our 
log--concavity assumption (\ref{logco}) 
we have $r_+(n) - r_+(n+1)\geq 0$, $n=0,1,\dots,N-1$. 
Moreover, by Lemma \ref{equi} $\tilde c(n+1)-\tilde c(n)\geq \d$ for some 
$\d>0$.
Therefore, by Lemma 
\ref{claim1} (with $\d_+=0$) we have the desired estimate (\ref{222})
with $\tilde c$ in place of $c$, and (\ref{222}) follows from the equivalence 
$\tilde c \asymp c$. \qed

\smallskip

Thanks to the equivalence 
$\nu_x\asymp\tilde\nu_x$ and (\ref{nutx}), the proof of
Proposition \ref{one-v} is an immediate consequence of Lemma  \ref{lele} if 
we can prove
\be
\mu_{L-1}\left(N-n\right)
\asymp \nep {-V(n)}\,,
\la{eqq}
\end{equation}
with a function $V$ satisfying (\ref{logco}). 
To prove (\ref{eqq}) we introduce the 
standard grand--canonical zero--range measures.
For every $\a>0$ and every vertex $x$ we consider the measures
\be
\mu_{x,\a}(0)
=\frac1{Z_\a}\,,\quad  \mu_{x,\a}(n)=
\frac{\a^n}{Z_\a}\prod_{k=1}^n 
\frac1{ c(k)}\,.
\la{mualfa}
\end{equation}
For every $\r>0$, let $\a_\r > 0$ denote the unique value of $\a$ such
that 
\be
\sum_{n=0}^\infty n\mu_{x,\a}(n)=\r\,.
\la{rhor}
\end{equation}
It is customary to write 
simply $\mu_{x,\r}$ for $\mu_{x,\a_\r}$. 
Similarly we denote by $\mu_{L,\r}$
the product $\otimes_{x\in V_L}\mu_{x,\r}$.
Setting $\r_n:=(N-n)/(L-1)$, for every $n\leq N-1$ we can write
\be
\mu_{L-1}\left(N-n\right)
= (\a_{\r_n})^{n-N} 
\, \left(\frac{Z_{\a_{\r_n}}}{Z_1}\right)^{L-1}\,
\mu_{L-1,\r_n}\left(N-n\right)\,.
\la{mul1}
\end{equation}
The idea is to use (\ref{mul1}) for all values of $n$ 
except those for which $N-n$ becomes too small. Therefore 
we fix an integer $m>0$, set $N_0=N-m$, and will use the identity
(\ref{mul1}) for all $n\leq N_0$. Here we proceed as follows.
Denoting by $\si^2_\r$ the variance of $\mu_{x,\r}$ we have the
following well known bounds 
see e.g.\ \cite{LSV,DP}:
\begin{gather}
\si^2_\r\asymp \r 
\la{asymp1}\\
\mu_{L,\r}\left( \r L\right) \asymp (\si_\r^2 L)^{-\frac12}
\,.
\la{asymp2}
\end{gather}
This implies $\mu_{L-1,\r_n}\left(N-n\right)\asymp (N-n)^{-\frac12}$.
Therefore from (\ref{mul1})
\be
\mu_{L-1}\left(N-n\right)
\asymp
\nep {- \tilde V(n)}\,,
\la{mul2}
\end{equation}
where, for every $t\in[0,N)$ we define $\r_t = (N-t)/(L-1)$ and  
\be
\tilde 
V(t)= (N-t)\log \a_{\r_t} - (L-1)\log (Z_{\a_{\r_t}}/Z_1) +\frac12 \log(N-t)\,.
\la{vt}
\end{equation}
We now prove that $\tilde V$ is convex if $t$ is not too close to $N$, 
i.e. $\tilde V''(t)\geq 0$, $t\leq N_0$. 
Clearly $\tilde V''(t) = \varphi''(t) - (2(N-t)^2)^{-1}$ with 
$\varphi(t):= (N-t)\log \a_{\r_t} - (L-1)\log (Z_{\a_{\r_t}}/Z_1)$.
We have $$\varphi'(t) = -\log(\a_{\r_t}) + (N-t)\dert\log(\a_{\r_t})
-(L-1)\dert\log(Z_{\a_{\r_t}})\,.$$
Using (\ref{rhor}) we see that 
$\dert\log(Z_{\a_{\r_t}}) = \r_t\dert\log(\a_{\r_t})$ and 
the last two terms in the expression for $\varphi'(t)$ cancel each other.
We then have $\varphi''(t) = - \dert\log(\a_{\r_t})$. Reasoning as above
and using $\dert \r_t = -1/(L-1)$ we have $\varphi''(t) = 1/(L-1)\si^2_{\r_t}$.
Therefore, for some
independent $C\in[1,\infty)$
\be
\tilde V''(t) 
= \frac{1}{(L-1)\si^2_{\r_t}} -\frac{1}{2(N-t)^2}
\geq \frac1{C(N-t)} -\frac{1}{2(N-t)^2}\,,
\la{20}
\end{equation}
where in the last estimate we have used (\ref{asymp1}).
Then $\tilde V''(t)\geq 0$ for all $N-t\geq C/2$. 
This implies -- by integration --
that $\grad \tilde V(n)\leq  \grad \tilde V(n+1)$
at least for all $n\leq N - 2 - C/2$. Setting e.g.\ $m = [C]$
we have shown that $\grad^2 \tilde V(n)\geq 0$, $n\leq N_0-2=N - m - 2$.

\smallskip

We still have to deal with the case $N-n\leq m$. 
Here we use the fact that 
\be
\mu_{L-1}(k) \asymp L^k \mu_x(0)^L\,,\quad \; k\leq m\,.
\la{mux0}
\end{equation}
To prove the lower bound in (\ref{mux0}) we simply 
observe that putting $k$ particles in $k$ different sites one has,
for some $m$--dependent $C<\infty$
$$
\mu_{L-1}(k)\geq \binom{L-1}{k} \mu_x(1)^k\mu_x(0)^{L-k}
\geq \frac1C\, L^k \mu_x(0)^L\,,\quad \; k\leq m\,.
$$
Similarly the upper bound is obtained by requiring at least $L-1-k$
sites to be empty:
$$
\mu_{L-1}(k)\leq \sum_{\ell\geq L-1-k} \binom{L-1}{\ell}
\mu_x(0)^{\ell} \leq C\,L^k \mu_x(0)^L\,,\quad \; k\leq m\,.
$$
Summarizing, from (\ref{mul2}) and (\ref{mux0}) 
we have obtained that, for every fixed $K\in(0,\infty)$
the equivalence (\ref{eqq}) holds 
with the function 
$V=V_K$ given by 
\be
V(n)=
\begin{cases}
\tilde V(n) & n\leq N_0\\
(n-N)\log L - L\log\mu_x(0) + K^{n-N_0} & N_0 < n \leq  N
\end{cases}
\la{vnn}
\end{equation}
Note that the addition of the term $K^{n-N_0}$ in
(\ref{vnn}) does not break the equivalence since $n-N_0\leq m$.
What we have seen in (\ref{20}) implies
$\grad^2V(n)\geq 0$ for $n\in [0,N_0 -2]$.
We are left with the case $n\geq N_0-1$. But this is
easily obtained by taking the constant $K$
sufficiently large. For instance: 
from (\ref{mux0}) we know that, for some universal constant $C<\infty$
$\tilde V(N_0) \geq (N_0 - N)\log L - L \log \mu_x(0) - C$,
so that $\grad V(N_0) \leq  \log L + K + C$. On the other hand
$\grad V(N_0+1)= \log L + K^{2} - K$. For $K$ large this gives 
$\grad^2V(N_0)\geq 0$. Similar reasoning applies for the remaining values 
of $n\geq N_0 - 1$. This ends the proof of the claim in (\ref{eqq}) and
concludes the proof of Proposition \ref{one-v}.

\section{Proof of Theorem \ref{mainre}}
The proof of Theorem \ref{mainre} is based on a
variant of the martingale recursive
method developed in \cite{LY}, see also \cite{LSV,LPY,GQ,LeeY,Yau,Djalil}. 
We set
\be
\g(L) = \sup_{N \geq 1} \g(L,N)\,.
\la{gel} 
\end{equation}
Note that the result of \cite{DP} on the logarithmic Sobolev
inequality for the local dynamics defined by (\ref{Dirich2})
implies that $\g(L)<\infty$ for every $L$. 
We are going to prove
\be
\la{i3}
\sup_L \,\g(L)<\infty\,.
\end{equation}
To this end we start with 
the usual decomposition of entropy and write, for $f>0$
\be
\la{i4}
\ent_{\nu}(f)=\frac1L\sum_{x}\nu\left[ \ent_{\nu}(f\tc \eta_x)\right]
+ \frac1L\sum_{x}\ent_{\nu}(f_x)\,,
\end{equation}
where $\ent_{\nu}(f\tc \eta_x)$ 
denotes the entropy of $f$ w.r.t.\ 
$\nu[\cdot\tc\eta_x]$ (the measure $\nu$ conditioned
to have a given number of particles $\eta_x$ at $x$)  
and we have defined
$$
f_x(\eta) = f_x(\eta_x) = \nu[f\tc\eta_x]\,.
$$
Since, for every given $0\leq \eta_x\leq N$, the measure
$\nu[\cdot\tc\eta_x]$ coincides with the canonical zero--range measure
on $L-1$ vertices with total particle number $N-\eta_x$,
we can estimate, for every $\eta_x$
$$
\ent_{\nu}(f\tc \eta_x) 
\leq \g(L-1)
\frac1{L-1}\sum_{y\neq x}\sum_{z\neq x}
\nu\left[c_y\,\grad_{yz}f\,\grad_{yz}\log f\tc \eta_x\right]\,.
$$
Taking $\nu$--expectation and 
averaging the above expression over $x$ we obtain that (\ref{i4})
is bounded above by
\be
\g(L-1)\,\frac{L-2}{L-1} \,\cE_{\nu}(f,\log f) + 
 \frac1L\sum_{x}\ent_{\nu}(f_x)\,.
\la{i40}  
\end{equation}
The next two subsections will explain how to estimate the second term in 
(\ref{i40}). Here we anticipate that the final result 
(see (\ref{one-v-bn}) and (\ref{bn-e}) below) 
will be that for every $\eps >0 $ there exist two
constants $\ell_\eps,C_\eps<\infty$ independent of $L$ and $N$ such that
for all $L\geq \ell_\eps$ we have
\be
\sum_{x}\ent_{\nu}(f_x) \leq \eps \,\ent_\nu(f) + C_\eps\,\cE_\nu(f,\log f)\,.
\la{claim}
\end{equation}
Once the above result is available it is easy to end the proof of
(\ref{i3}). Indeed, from (\ref{claim}) and (\ref{i40})
we obtain
\be
\left(1-\frac\eps{L}\right)\,\g(L)\, \leq\, \frac{L-2}{L-1} \,\g(L-1) + 
\frac{C_\eps}{L} \,,\quad\, L\geq \ell_\eps\,,\la{i80}
\end{equation}  
which implies the claim (\ref{i3}) if $\eps$ is sufficiently small
(e.g.\ $\eps<\frac12$).

%

\subsection{From one--vertex estimate to covariances}
Let us recall the following change of variable relation:
for any function $f$ and any pair of vertices $x,y$ 
\be
\nu[c_x f] = \nu[c_y f^{yx}] \,.
\la{m2}
\end{equation}
The above is an immediate consequence
of the definitions of the symbols involved and 
the fact that, for any $\eta\in\O_L$ with $\eta_x\geq 1$ we have 
$$
\frac{\nu(\eta^{xy})}{\nu(\eta)} = 
\frac{\mu_x(\eta_x-1)\mu_y(\eta_y+1)}{\mu_x(\eta_x)\mu_y(\eta_y)} = 
\frac{c_x(\eta_x)}{c_y(\eta_y +1)}\,.
$$

We start our proof of the claim (\ref{claim}) with an application of Proposition \ref{one-v}
to the function  $u=f_x/\nu[f_x]$. Here and in the rest of this subsection $x$
is an arbitrary fixed vertex. 
We have 
\be
\ent_{\nu}(f_x) \leq 
C\, \sum_{n=0}^N
\nu_x(n) c(n) [f_x(n)- f_x(n-1)] \log \frac{f_x(n)}{f_x(n-1)}\,.
\la{one-v-f}
\end{equation}
To estimate the R.H.S.\ of (\ref{one-v-f}) we first rewrite things as follows.
For every vertex $y\neq x$ and
for every $n$ we define the functions
\be
g_{x,y,n}(\eta)=\frac{c_y(\eta)}{\nu[c_y\tc \eta_x=n]}\,,
\quad g_{x,n}(\eta) = \frac{1}{L-1}\sum_{y\neq x}g_{x,y,n}(\eta)\,.
\la{m00}
\end{equation}
In order to simplify notations, below we will 
write $\nu[\cdot\tc n]$ for $\nu[\cdot\tc \eta_x=n]$. 
Formula (\ref{m2}) can be used to deduce
the identity 
\be 
f_x(n)=\nu[f\tc n]=\nu[g_{x,y,n-1} f^{yx}\tc n-1]\,,
\la{m11} 
\end{equation}
valid for every $y\neq x$ and $n\geq 1$. Indeed, 
write $\chi_{x,n}(\eta)$ for the indicator function of
the event $\{\eta_x=n\}$. Then $(\chi_{x,n})^{yx} = \chi_{x,n-1}$
and 
$$
\nu[f\chi_{x,n}] =
\frac1{c(n)}\, \nu[c_x f\chi_{x,n}] = 
\frac1{c(n)}\,\nu[c_y f^{yx}\chi_{x,n-1}]\,.
$$
When $f=1$ this shows that
$\nu[c_y\tc n-1]=\frac{c(n)\nu_x(n)}{\nu_x(n-1)}$ and (\ref{m11}) 
follows. 

In particular, (\ref{m11}) shows that
\be
\la{m111}
\nu[f\tc n] - \nu[g_{x,n-1}f\tc n-1] = 
\frac{1}{L-1}\sum_{y\neq x}
\nu\left[g_{x,y,n-1}\grad_{yx}f \tc n-1\right]\,.
\end{equation}
Our first step in the estimate of (\ref{one-v-f}) is the next lemma. 
We recall the standard notation $\mu[f,g]=\mu[fg]-\mu[f]\mu[g]$ for the covariance
of two functions $f,g$ w.r.t.\ a measure $\mu$. 
\begin{Le} 
\la{prom} 
There exists $C<\infty$ such that for every $f>0$,  $L\geq 2$, and $N\geq n\geq 1$
\be 
[f_x(n)-f_x(n-1)]\log \frac{f_x(n)}{f_x(n-1)}\leq C\, \big\{A_x(n) \,+\,  B_x(n)\big\}\,,
\la{anbn}
\end{equation}  
where we define
\begin{gather}
A_x(n)= \left(\nu[f\tc n]-\nu[g_{x,n-1} f\tc n-1]\right)\,\log\frac{\nu[f\tc n]} 
{\nu[g_{x,n-1} f \tc n-1]}\,,
\la{an}\\
B_x(n)= \frac{\nu[g_{x,n-1},f\tc n-1]^2}
{f_x(n)\vee f_x(n-1)}\,.
\la{bn}
\end{gather}
\end{Le} 
\proof
Set $a=f_x(n)$, 
$b=\nu[g_{x,n-1} f\tc n-1]$ and $c=f_x(n-1)$.
With the notation $\a(a,b)=(a-b)\log(a/b)$, the desired estimate (\ref{anbn})
can be written as 
\be
\a(a,c) \leq C\,\a(a,b)
+ C\,\frac{(b-c)^2}{a\vee c}\,.
\la{m4}
\end{equation}
Note that the above inequality cannot hold for all $a,b,c>0$
without restrictions (take e.g.\ $c=1$, $a=b$ and let $b\nearrow\infty$). 
The point is that in our setting we
have $1/C\leq b/c \leq C$, for some possibly different $C\in[1,\infty)$. To see this
recall that $\nu[\eta_y\tc n]=(N-n)/(L-1)$ for all $n$ and $y\neq x$ and 
use $c(n)\asymp n$ to obtain
\be
g_{x,n-1}(\eta)\asymp\frac{\sum_{y\neq x}\eta_y}{N-(n-1)}\,= 1 \,,\quad
\nu[\cdot\tc n-1]\,-\,a.s.\
\la{m5}
\end{equation}
for every $n\geq 1$. 
Therefore $b\asymp c$.
We now write
\begin{gather*}
\a(a,c)=c(a/c -1)\log(a/c) = c \,h(t)\,,\\
\quad\; h(t):=t(\nep{t}-1)\,,
\quad t:=\log(a/c)\,.
\end{gather*}
It is not difficult to check the function $h$ satisfies: for every $C<\infty$ 
\begin{gather}
  \sup_{u\leq C}\,\sup_{t\leq 2C}\frac{h(t)}{h(t-u)+u^2}\,<\,\infty\,, 
\la{lemh1} \\
  \sup_{u\leq C}\,\sup_{t\geq 2C}
  \,\frac{h(t)}{h(t-u)}\,<\,\infty\,. \la{lemh2}
\end{gather}

In the rest of this proof we use $C_1,C_2,\dots$ to denote finite positive 
constants (independent of $n,N,L$). 
Setting $u:=\log(b/c)$, we 
know that $u \leq C_1$. Suppose first that $a/c\leq 2C_1$.
Then by  (\ref{lemh1})
we know that there exists $C_2<\infty$
such that $h(t)\leq C_2(h(t-u) + u^2)$, i.e.
$$
\a(a,c)\leq C_2\,\left[c(a/b -1)\log(a/b) +c (\log(b/c))^2\right]\,.
$$
The first term above is $c/b\,\a(a,b)\leq C \a(a,b)$. 
For the second term we use the elementary fact that for every 
$\d > 0$, there is $C=C(\d)<\infty$ such that 
$|\log(1+x)|\leq C\, |x|$, for any $x\geq \d - 1$. 
With $x=b/c -1$, this says that the second term 
is bounded by 
$$
C_3\,\frac{(b-c)^2}{c} \,\leq C_4\, \frac{(b-c)^2}{a\vee c}\,,
$$
where we used the assumption $a/c \leq 2C_1$.
This completes the proof of (\ref{m4}) under this assumption. 
If $a/c > 2C_1$ we have by  (\ref{lemh2}) 
$h(t)\leq C_5 h(t-u)$, i.e.\
$$
\a(a,c)\leq C_5 \, c\, (a/b -1)\log(a/b) \leq C_6 \a(a,b)\,,
$$
which clearly implies (\ref{m4}). \qed
\smallskip


When we insert the estimate of
Lemma \ref{prom} in (\ref{one-v-f}) we therefore obtain two terms, corresponding
to $A_x(n)$ and $B_x(n)$, respectively. We explain here how to bound the first term. 
This is a modification of a
rather standard convexity argument, see e.g.\ \cite{GQ}.
The more delicate estimate of the term coming from $B_x(n)$ is given in the next subsection.

\smallskip

Thanks to the identity (\ref{m11}) and the convexity of
$(a,b)\to (a-b)\log(a/b)$ on $\bbR_+\times \bbR_+$, Jensen's inequality implies
\be
A_x(n)\leq \frac{1}{L-1}\sum_{y\neq x}
\nu\left[g_{x,y,n-1}\grad_{yx}f\,\grad_{yx}\log f\,\tc n-1\right]\,.
\la{m20}
\end{equation}
Going back to (\ref{one-v-f}) and
using (see (\ref{m11})) 
$$
\nu_x(n)c(n)g_{x,y,n-1} =  \nu_x(n-1)\,c_y\,,\quad y\neq x\,,
$$ 
we see that 
\be
\sum_{n=0}^N
\nu_x(n) c(n) A_x(n)\leq 
\frac1{L-1}\sum_{y\neq x} \nu\left[c_y\,\grad_{yx}f\,\grad_{yx}\log f\right]\,.
\la{m30}
\end{equation}
When we sum over $x$ in (\ref{one-v-f}), from
Lemma \ref{prom} and (\ref{m30}) we obtain
\be
\sum_x\ent_{\nu}(f_x) \leq 
C\, \cE_\nu(f,\log f) + 
C\,\sum_x\sum_{n=0}^N
\nu_x(n) c(n) B_x(n) \,.
\la{one-v-bn}
\end{equation}

\subsection{The covariance estimate}
We need the following key estimate
on covariances. 
\begin{Pro}
\la{procova}
Assume (H.1) and (H.2). For  every $\eps>0$, there exist finite constants $C_\eps$ and $\ell_\eps$
such that for every $L\geq \ell_\eps$, $N\geq 1$ and for every $f>0$
\be
\nu\left[f,\sum_x c_x\right]^2 \leq N\,\nu[f]\,\left[C_\eps \cE_\nu(\sqrt f,\sqrt f) + \eps \,\ent_\nu(f)
\right]\,.
\la{procova1}
\end{equation}
\end{Pro}

Before going to the proof we want to make sure this
result is indeed sufficient for our claim (\ref{claim}) to hold.
To this end we fix a vertex $x$ and  
apply (\ref{procova1}) by replacing $\nu$ with $\nu[\cdot \tc n-1]=\nu[\cdot \tc \eta_x=n-1]$,
$L$ by $L-1$ and $N$ by $N-n+1$. 
Using the equivalence 
\be
\nu[c_y\tc n-1] \asymp \frac{N-n+1}{L-1}\,,
\la{equi10}
\end{equation}
we then see that for some $C<\infty$
\be 
\nu[f,g_{x,n-1}\tc n-1]^2 \leq  \frac{C\,f_x(n-1)}{N-n+1}
\left[C_\eps \cE_\nu(\sqrt f,\sqrt f\tc n-1) + \eps \,\ent_\nu(f\tc n-1)
\right]\,.
\la{procovo}
\end{equation}
Using again (\ref{equi10}) and the identity $\nu_x(n)c(n)=\nu[c_y\tc n-1]\nu_x(n-1)$
we get, with a possibly different constant $C$
\be
\sum_{n=0}^N
\nu_x(n) c(n) B_x(n) 
\leq \frac{C}L \,
\left[C_\eps \cE_\nu(\sqrt f,\sqrt f) + \eps \,\ent_\nu(f)
\right]\,,
\la{bn-e}
\end{equation}
where we have used the easily verified
estimates $$\nu[\cE_\nu(\sqrt f,\sqrt f\tc \eta_x)]\leq 
\cE_\nu(\sqrt f,\sqrt f)\,,\quad 
\nu[\ent_\nu(f\tc \eta_x)]\leq 
\ent_\nu(f)\,.
$$ 
Finally, the desired estimate (\ref{claim}) follows from
(\ref{bn-e}), (\ref{one-v-bn}) and the elementary bound (\ref{deco}).

\smallskip

We turn to the proof of Proposition \ref{procova}. 
%
Let us first recall the covariance estimate proved in \cite{DP}.
Corollary~3.11 there states that assuming (H.1) and (H.2) one has 
\be
\nu\left[f,\sum_x c_x\right]^2 
\leq C\,N\,\nu[f]\,\left[\nu[f]+ C_\eps \,L^2\,\cD_\nu(\sqrt f,\sqrt f) + \eps \ent_\nu(f)
\right]\,.
\la{procova2}
\end{equation}
Here $C$ is a finite constant depending only on the
parameters appearing in (H.1) and (H.2) and 
$\cD_\nu$ stands for the local Dirichlet form defined in (\ref{Dirich2}).
The constants $\eps$ and $C_\eps$ have the same meaning as in our Proposition \ref{procova} above. 
To prove our bound in (\ref{procova1}) we therefore 
have to improve the latter result in two ways:
first, we need to replace $L^2\,\cD_\nu(\sqrt f,\sqrt f)$ by $\cE_\nu(\sqrt f,\sqrt f)$
and second, we have to remove the extra term $\nu[f]$ appearing in (\ref{procova2}).
It turns out that the first improvement requires only straightforward modifications 
of the argument of \cite{DP}. The second, on the other hand, will require 
some additional work, which will be based on a combination of ideas from \cite{DP} and \cite{LPY}.
As in \cite{DP} we consider separately the case of small density and
the case of densities uniformly bounded away from zero. 
In the rest of the proof of 
Proposition \ref{procova} we adopt the convention that $C$ represents 
a generic finite constant which may only depend on the 
parameters appearing in (H.1) and (H.2). When constants depend  
on a further parameter as e.g.\ $\eps,M$ or $K$ 
we write this explicitly as $C_\eps,C_M$ or $C_K$ respectively. 
In all cases it is understood that these constants are independent of $L$ and $N$.
We warn the reader that the numerical value of these constants may change from
line to line.
%
\subsection{Small density}
Here we assume that $\r:=N/L$ satisfies $\r\leq \r_0$ with $\r_0$ 
a parameter to be taken sufficiently small depending on $\eps$.
Recall the definition (\ref{rhor}) of the parameter $\a_\r$.
We use the notations 
\begin{gather}
\varphi_x(\eta_x)= c(\eta_x) - \frac{\a_\r}{\r}\,\eta_x\nonumber\\
\bar \varphi_x = \varphi_x-\nu[\varphi_x]\,, 
\quad\Phi(\eta) = \sum_x\bar \varphi_x\,.
\la{varphi}
\end{gather}
\begin{Le}
\la{exptfi}
For every $M>0$, there exists $C_M<\infty$ such that 
\be
\frac1t \log\nu\left[\,\exp t|\Phi|\,\right]\leq C_M N\sqrt\r \,t\,,\quad \; t\in[0,M]\,.
\la{exp1}
\end{equation}
\end{Le}
Before giving a proof we show that Lemma \ref{exptfi} implies that for
every $M>0$, there exists $C_M<\infty$ such that for any $f>0$, with $\nu[f]=1$:
\be
\nu\left[f,\sum_x c_x\right]^2 
\leq C\,N\,\left[(C_M \r) \vee \frac1M \right]
\ent_\nu(f)\,.
\la{prosmall}
\end{equation}
Of course, 
by taking $\r_0$ small enough, (\ref{prosmall}) gives the desired result 
(\ref{procova1}) for small density.
To prove (\ref{prosmall}) we use the entropy inequality to write, for every $t>0$
$$
\nu\left[f,\sum_x c_x\right] = \nu\left[f\Phi\right]
\leq \frac1t \,\log\nu\left[\,\exp t \Phi\,\right] + \frac1t\,\ent_\nu(f)\,.
$$
We may apply the above inequality with $-\Phi$ replacing $\Phi$.  
Therefore, passing to absolute values, Lemma \ref{exp1} gives
\be
\left|\nu\left[f,\sum_x c_x\right]\right| \leq C_M N\sqrt\r \,t
+ \frac1t\,\ent_\nu(f) \,,\quad \; t\in[0,M]\,.
\la{prosmall1}
\end{equation}
Set now $\bar t = \sqrt{\frac{M}{N}\,\ent_\nu(f)}$. 
If $\bar t \leq M$, (\ref{prosmall}) follows 
immediately by plugging $t=\bar t$ in (\ref{prosmall1}). 
If, however, $\bar t \geq M$ we may use the rough bound $|\sum_x c_x|\leq C\,N$ 
to estimate
$$
\nu\left[f,\sum_x c_x\right]^2\leq C\,N^2\leq C N\,\frac1M\,\ent_\nu(f)\,.
$$
We now turn to the proof of Lemma \ref{exptfi}. Since all our estimates
below are easily seen to hold with $\Phi$ replaced by $-\Phi$ we may 
restrict to estimate $\nu\left[\exp t\Phi\right]$ instead of 
$\nu\left[\exp t|\Phi|\right]$. We consider two different cases: $M\geq t\geq M\wedge\sqrt L /N$ 
and $t\leq M\wedge 
\sqrt L /N$.

\smallskip

\noindent
{\bf Case $M\geq t\geq M\wedge(\sqrt L /N)$}. We recall the following bound derived in \cite{DP},
see (4.80) and (4.88) there:
\be
\frac1t \log\nu\left[\,\exp t\Phi\,\right]\leq 
\frac{C}t\,+\, C\sqrt N\,+\,C_M N\,\r \,t \,,\;\quad \; t\in[0,M].
\la{exp3}
\end{equation}
If $t\geq \sqrt L /N$ we have $1/t \leq  N \r \,t$ and $\sqrt N\leq N\sqrt \r\, t$.
Therefore (\ref{exp1}) is contained in (\ref{exp3}) in this case.

\smallskip

\noindent
{\bf Case $t\leq M\wedge (\sqrt L /N)$}. 
The bound (\ref{exp3}) is not optimal for small values of $t$
and we need a different approach here. We may proceed as in \cite{LPY}, Lemma 6.5. 
Without loss of generality, we assume that $L$ is even. We call $V_{L/2}$
the set of vertices $\{1,2,\dots,L/2\}$.  
By Schwarz inequality we
have
$$
\log \nu\left[\exp t\Phi\right] \leq \log \nu\left[\exp 2t\widetilde\Phi\right]\,,
\quad \widetilde\Phi:=\sum_{x\in V_{L/2}}\bar\varphi_x\,.
$$
For every function $g$ such that $\nu[g]=0$ we may estimate
%
\be
\nu[\nep{g}]\leq \exp\left\{
\frac12\,\nu\left[g^2\nep{|g|}\right]\right\}\,.
\la{expexp}
\end{equation}
This estimate follows from $\nep a\leq 1+a+\frac12 a^2\nep{|a|}$, and
$1+x\leq \nep x$. 
We apply this bound to $g=2t\widetilde\Phi$.
Using the equivalence of ensembles bound (see e.g.\ Proposition~4.1 in \cite{DP})
we have $\nu[\widetilde\Phi^2\exp 2t|\widetilde\Phi|]\leq 
C\mu_{L,\r}[\widetilde\Phi^2\exp 2t|\widetilde\Phi|]$ 
and therefore 
\be
\nu\left[\exp 2t\widetilde\Phi\right]\leq 
\exp\left\{
C\,t^2\,\mu_{L,\r}\left[\widetilde\Phi^2\nep{2t|\widetilde\Phi|}\right]\right\}  \,.
\la{exp4}
\end{equation}
All the estimates below can be obtained for $-\widetilde\Phi$ as well as for $\widetilde\Phi$
without any change, 
therefore we will restrict to bound the expression
\be
\mu_{L,\r}\left[\widetilde\Phi^2\nep{2t\widetilde\Phi}\right]
= \sum_{x,y\in V_{L/2}} \mu_{L,\r}\left[\bar \varphi_x
\bar \varphi_y \nep{2t\widetilde\Phi}\right] = E_1 + E_2\,,
\la{exp5}
\end{equation}
where, using the product structure of $\mu_{L,\r}$ and writing $\mu_\r=\mu_{1,\r}$
$$
E_1 : = \frac{L}2\,
\mu_\r\left[\bar \varphi_1^2\nep{2t\bar\varphi_1}\right]
\mu_\r\left[\nep{2t\bar\varphi_1}\right]^{\frac{L}2-1}\,,
$$
$$
E_2 : = \frac{L}2\left(\frac{L}2-1\right)\,
\mu_\r\left[\bar \varphi_1\nep{2t\bar\varphi_1}\right]^2
\mu_\r\left[\nep{2t\bar\varphi_1}\right]^{\frac{L}2-2}\,.
$$
Recalling (see e.g.Corollary~6.4 in \cite{LSV}) that $|\bar\varphi_1 - \varphi_1|\leq C\frac{\sqrt{1+\r}}{L}$,
we estimate $$\mu_\r\left[\nep{2t\bar\varphi_1}\right]\leq \nep{t \frac{C}{L}}
\mu_\r\left[\nep{2t\varphi_1}\right]\,.$$ 
From (4.88) in \cite{DP}, 
$\mu_\r\left[\nep{2t\varphi_1}\right]\leq \nep{C_M\r^2 t^2}$, $t\leq M$. Therefore
$$
\mu_\r\left[\nep{2t\bar\varphi_1}\right]^{\frac{L}2-1}\leq C_M\, \nep {C_M\r^2t^2 L}
\leq C_M\,,
$$
the last bound following from $t^2\leq L/N^2$. 
This gives $E_1\leq C_M L \,\mu_\r
\left[\bar \varphi_1^2\nep{2t\bar\varphi_1}\right]$. Replacing as above $\bar\varphi_1$ with
$\varphi_1$ we have
\be
\mu_\r
\left[\bar \varphi_1^2\nep{2t\bar\varphi_1}\right]
\leq \frac{C}{L^2} + C\, \mu_\r\left[ \varphi_1^2\nep{2t\varphi_1}\right]\,.
\la{exp6}
\end{equation}
By direct computation (or reasoning as in (4.82),(4.84) and (4.86) in \cite{DP})
it is not hard to obtain the bound
\be
\mu_\r\left[ \varphi_1^2\nep{2t\varphi_1}\right]\leq C_M\r^2\,,\quad t\leq M\,.
\la{exp7}
\end{equation}
From (\ref{exp6}) and (\ref{exp7}), using $\r\geq 1/L$, we have obtained $E_1\leq C_M N\,\r$.
We now look for a similar bound on $E_2$. We first observe that for any $a\in\bbR$
we have $a\nep a\leq a + a^2\nep{|a|}$. Setting $a = 2t\bar\varphi_1$ we
obtain
$$
\mu_\r\left[\bar \varphi_1\nep{2t\bar\varphi_1}\right]
\leq \mu_\r\left[\bar \varphi_1\right] + 2t
\mu_\r\left[\bar \varphi_1^2\nep{2t|\bar\varphi_1|}\right]\,.
$$
Estimating as in (\ref{exp6}) and (\ref{exp7}) once for $\bar\varphi_1$
and once for $-\bar\varphi_1$, the second term above is bounded by $C_Mt\,\r^2$.
Since $\mu_\r[\varphi_1]=0$, direct computations show that
$|\mu_\r[\bar\varphi_1]|\leq C(\r^2\wedge \frac1L)$. Therefore
$$\mu_\r\left[\bar \varphi_1\nep{2t\bar\varphi_1}\right]
\leq C \frac\r{\sqrt L} + C_Mt\,\r^2\,.$$ 
Reasoning as above it is not hard to check that the last 
estimate holds for $-\mu_\r\left[\bar \varphi_1\nep{2t\bar\varphi_1}\right]$ as well.
We then obtain 
$$\mu_\r\left[\bar \varphi_1\nep{2t\bar\varphi_1}\right]^2
\leq C \frac{\r^2}{L}  + C_Mt^2\,\r^4\,.$$
This implies the estimate $E_2\leq C \,L\,\r^2 + C_M L^2\,t^2\,\r^4\,.$
Using the constraint $t^2\leq L/N^2$ this becomes 
$E_2\leq C_M N\,\r\,.$ In conclusion: from (\ref{exp4}) and (\ref{exp5})
we have
\be
\frac1t \log\nu\left[\,\exp t\Phi\,\right]\leq 
C_M N\,\r \,t\,. 
\la{exp9}
\end{equation}
This ends the proof of Lemma \ref{exptfi}.
%
%
%
\subsection{Density bounded away from zero} 
To prove Proposition \ref{procova} in the regime $\r\geq \r_0$
we need the following standard coarse graining procedure. We fix a parameter $K>0$
to be taken sufficiently large in the sequel. Without loss of generality we will assume that 
$K$ divides $L$ so that the set of vertices $V_L$ is the disjoint union of $\ell:=L/K$ 
sets of vertices $B_1,\dots,B_\ell$, each of cardinality $K$. 
We write $N_j = N_j(\eta)=\sum_{x\in B_j}\eta_x$ for the number of particles in the block $B_j$
and write $\cG$ for the $\si$--algebra generated by the functions $\eta\to N_j(\eta)$, $j=1,\dots,\ell$. 
In this way, the conditional expectation $\nu[\cdot\tc \cG]$ becomes the product 
$\prod_{j=1}^\ell\nu_{j,N_j}[\cdot]$,
where $\nu_{j,N_j}$ denotes the canonical zero--range measure on the $j$--th block
with $N_j$ particles. We start with the decomposition
\be
\nu\left[f,\sum_x c_x\right] = \nu\left[\nu\Big[f,\sum_x c_x\tc \cG\Big]\right]
+ \nu\left[f,\sum_{j=1}^\ell\nu_{j,N_j}
\Big[\sum_{x\in B_j} c_x\Big]\right]\,.
\la{covax}
\end{equation}
As in \cite{DP}, Corollary 3.9, it is not hard to prove 
\be
\nu\left[\nu\Big[f,\sum_x c_x\tc \cG\Big]\right]\leq 
C_KN \,\nu[f]\,\cE_\nu(\sqrt f,\sqrt f)\,.
\la{rough}
\end{equation}
We now concentrate on a bound on the second term in (\ref{covax}).
To this end we introduce the following notations. For every $x$ we set
$\bar c_x(\eta)=c(\eta_x) - \a'_\r\eta_x$, where $\a'_\r=\frac\dd{\dd\r}\a_\r$, and,
with $\r_j:=N_j/K$, for ever $x\in B_j$
\begin{gather*}
G(\r_j)=\nu_{j,N_j}[\bar c_x] - \mu_{x,\r}[\bar c_x]\,,\\
\bar G(\r_j)=\nu_{j,N_j}[\bar c_x] - \nu[\bar c_x]\,.
\end{gather*}
Note that these definition do not depend on the chosen $x\in B_j$. Moreover, 
$\mu[G(\r_j)]=0$ and $\nu[\bar G(\r_j)]=0$. 
We also set 
$$
\Psi(\eta)=K\sum_{j=1}^\ell \bar G(\r_j)\,,
$$
so that the second term in (\ref{covax}) becomes $\nu[f\Psi]$. 
Therefore our ultimate claim now becomes: for every $\r_0>0$, for every $\eps>0$ 
there exist constants
$K_\eps,\ell_\eps,C_\eps<\infty$ such that for all $K\geq K_\eps$
$\ell\geq \ell_\eps$
\be
\nu[f\Psi]^2\leq C\,N\,\nu[f]\,\left(C_\eps\cE_\nu(\sqrt f,\sqrt f) + \eps\,\ent_\nu(f)\right)\,.
\la{proso}
 \end{equation}
As a simplifying rule we do not write explicitly the $\r_0$--dependence of the various
constants.
\begin{Le}
\la{exptpsi}
For every $M>0$ there exists $C_{M}<\infty$ such that 
\be
\frac1t \log\nu\left[\,\exp t|\Psi|\,\right]\leq \frac{C_M N}{\sqrt K} \,t\,,\quad \; 0\leq t\leq 
\frac{M}{K\sqrt \r}\,.
\la{expa1}
\end{equation}
\end{Le}
Before giving the proof of the lemma we want to show that 
the estimate (\ref{expa1}) is sufficient to prove (\ref{proso}). 
As in (\ref{prosmall1}), assuming 
$\nu[f]=1$, (\ref{expa1}) allows to estimate
\be
\left|\nu\left[f\Psi\right]\right| \leq \frac{C_M N}{\sqrt K} \,t 
+ \frac1t\,\ent_\nu(f) \,,\quad \; 0\leq t\leq 
\frac{M}{K\sqrt \r}\,.
\la{pros1}
\end{equation}
Set again $\bar t = \sqrt{\frac{M}{N}\,\ent_\nu(f)}$. 
If $\bar t \leq \frac{M}{K\sqrt \r}$,  
plugging $t=\bar t$ in (\ref{pros1}) we have
\be 
\nu[f\Psi]^2 \leq C\,N\,\left(\frac{C_M}{K}\vee \frac{1}{M}\right)\,\ent_\nu(f) \,.
\la{pros3}
\end{equation}
Taking $M$ and $K$ sufficiently large in a suitable way this clearly
implies (\ref{proso}). The case $\bar t\geq \frac{M}{K\sqrt \r}$ is much more
delicate. By repeating exactly the computations in \cite{DP}, see (4.76) there,
in this case 
one arrives at the desired estimate (\ref{proso}) 
except that $\cE_\nu(\sqrt f,\sqrt f)$ is replaced by $L^2\cD_\nu(\sqrt f,\sqrt f)$.
To see how this can be improved we recall that the relevant term
comes from expressions (4.55) and (4.69) in \cite{DP}. In particular,
now the precise estimate we need in order to obtain our claim can be written as
\be
\left(\frac1\ell\sum_{i,j=1}^\ell\nu\left[\sqrt{(N_i+N_j)\,\nu_{i,j}[f]\,\cE_{i,j}(\sqrt f,\sqrt f)}\right]
\right)^2
\leq C_K N\,\cE_\nu(\sqrt f,\sqrt f)\,,
\la{pros4}
\end{equation}
where, following \cite{DP}, we write $\nu_{i,j}[f]=\nu[f\tc \cF_{i,j}]$, with $\cF_{i,j}$
denoting the $\si$--algebra generated by $\{\eta_x, x\in(B_i\cup B_j)^c\}$. 
Here $\cE_{i,j}$ stands for the 
Dirichlet form 
$$
\cE_{i,j}(\sqrt f,\sqrt f) = \sum_{x,y\in B_i\cup B_j}\nu_{i,j}[c_x(\grad_{xy}\sqrt f)^2]\,.
$$
To prove (\ref{pros4}) we observe that by Schwarz inequality
for the combined measure $\frac1{\ell^2}\sum_{i,j}\nu[\cdot]$, the L.H.S.\ of (\ref{pros4})
is bounded by 
$$
\left(\frac1{\ell}\sum_{i,j=1}^\ell\nu\Big[(N_i+N_j)\,\nu_{i,j}[f]\Big]\right)
\left(\frac1{\ell}\sum_{i,j=1}^\ell\nu\left[\cE_{i,j}(\sqrt f,\sqrt f)\right]\right)
$$
The first term above is handled by observing that $\nu\left[(N_i+N_j)\nu_{i,j}[f]\right]=
\nu\left[(N_i+N_j)f\right]$ and $\sum_{i}^\ell\nu\left[N_i f\right]=N$, since 
we are assuming $\nu[f]=1$.
Therefore (\ref{pros4}) follows from the following estimate, which is easily verified
$$
\frac1{\ell}\sum_{i,j=1}^\ell\nu\left[\cE_{i,j}(\sqrt f,\sqrt f)\right]\leq C_K\cE_\nu(\sqrt f,\sqrt f)\,.
$$
This ends the proof of (\ref{proso}) assuming the result of Lemma \ref{exptpsi}.

\smallskip

\noindent
{\bf Proof of Lemma \ref{exptpsi}}.
As in the proof of Lemma \ref{exptfi} we need to consider two regimes for the
values of $t$. 

\smallskip

\noindent
{\bf Case $\frac{M}{K\sqrt \r}\geq t\geq  \sqrt\frac{K}{N}$}. 
We use the following bound derived in \cite{DP},
see (4.33) there:
\be
\frac1t \log\nu\left[\,\exp t\Psi\,\right]\leq 
\frac{C}t\,+\, C\sqrt N\,+\, \frac{C_MN}{K} \,t \,,\;\quad \; 0\leq t\leq\frac{M}{K\sqrt \r}\,.
\la{expa5}
\end{equation}
If $t\geq \sqrt\frac{K}{N}$ we have $1/t \leq  \frac{N}{K}\,t$ and 
$\sqrt N\leq \frac{N}{\sqrt K}\, t$,
therefore (\ref{expa1}) is contained in (\ref{expa5}) in this case.

\smallskip

\noindent
{\bf Case $t\leq \frac{M}{K\sqrt \r}\wedge\sqrt\frac{K}{N}$}. 
In this case we use the same strategy as in Lemma \ref{exptfi}, in the case 
of small $t$.
The function $\Psi$ replaces now the function $\Phi$,
and the functions $K\bar G(\r_j)$ play here the role of the functions $\bar \varphi_x$ defined
in (\ref{varphi}).                                                        
We only sketch the arguments required to prove the needed estimates
since they are essentially the same as in the case of small density.
As in that case we may reduce the proof to suitable bounds on the expressions
$$
E_1 : = \frac{\ell}2\,
\mu_{K,\r}\left[(K\bar G(\r_1))^2\nep{2tK\bar G(\r_1)}\right]
\mu_{K,\r}\left[\nep{2tK\bar G(\r_1)}\right]^{\frac{\ell}2-1}\,,
$$
$$
E_2 : = \frac{\ell}2\left(\frac{\ell}2-1\right)\,
\mu_{K,\r}\left[K\bar G(\r_1)\nep{2tK\bar G(\r_1)}\right]^2
\mu_{K,\r}\left[\nep{2tK\bar G(\r_1)}\right]^{\frac{\ell}2-2}\,.
$$
We recall that (see e.g.\ Corollary~6.4 in \cite{LSV})
\be
|\bar G(\r_1)- G(\r_1)| = |\nu[c_x]-\mu_{x,\r}[c_x]| \leq C\frac{\sqrt{1+\r}}{L}\,.
\la{expa10}
\end{equation} 
Therefore, using $\r\geq \r_0$
\be
\mu_{K,\r}\left[\nep{2tK\bar G(\r_1)}\right]
\leq \nep{C K \sqrt \r t /L}\,\mu_{K,\r}\left[\nep{2tK G(\r_1)}\right]\,.
\la{expa11}
\end{equation} 
Moreover as in (\ref{expexp})
\be
\mu_{K,\r}\left[\nep{2tK G(\r_1)}\right] \leq \exp \left\{
2t^2 \mu_{K,\r}\left[(K G(\r_1))^2\nep{2tK |G(\r_1)|}\right]
\right\}\,.
\la{expa120}
\end{equation} 
An adaptation of estimates (4.12), (4.19) and (4.27) in \cite{DP}
yields the following crucial bound:
\be
\mu_{K,\r}\left[(K G(\r_1))^2\nep{2tK |G(\r_1)|}\right]
\leq C_M\,\r\,,\;\quad \; t\leq \frac{M}{K\sqrt \r}\,.
\la{expa12}
\end{equation} 
Since $t\leq\sqrt\frac{K}{N}$, (\ref{expa11}), (\ref{expa120}) and (\ref{expa12}) give
\be
\mu_{K,\r}\left[\nep{2tK\bar G(\r_1)}\right]^{\frac{\ell}2-1}\leq C_M\,.
\la{expa13}
\end{equation} 
Using again (\ref{expa10}) we see that (\ref{expa11}) and 
(\ref{expa12}) imply
\be
\mu_{K,\r}\left[(K\bar G(\r_1))^2\nep{2tK\bar G(\r_1)}\right]
\leq C_M\,\r\,.
\la{expa14}
\end{equation} 
Summarizing, we have obtained $E_1\leq C_M\,\ell\,\r = C_M\,N/K$. 
The estimate on $E_2$ can be done in the same way as we did for the case
of small density. In particular, using (\ref{expa14}) we obtain
$E_2\leq C_M\ell^2(\,K^2\,\r/L^2 + \,t^2\,\r^2)$. Since $t^2\leq K/N$
this gives $E_2\leq C_M \,N/K$. Therefore 
$$
\frac1t \log\nu\left[\,\exp t|\Psi|\,\right]\leq 
\frac{C_M N}{K} \,t\,,\quad \; 0\leq t\leq \frac{M}{K\sqrt \r}\wedge\sqrt\frac{K}{N}\,.
$$
This ends the proof of Lemma \ref{exptpsi}.
\qed

\smallskip
\bigskip

\noindent
{\bf Acknowledgments.}
We thank Paolo Dai Pra for several useful comments and discussions.
We acknowledge the support of M.I.U.R.(Cofin).

\end{document}